\documentclass{article}

\usepackage[english]{babel}
\usepackage{comment}
\usepackage[letterpaper,top=2cm,bottom=2cm,left=3cm,right=3cm,marginparwidth=1.75cm]{geometry}

\usepackage{amssymb, amsmath, amsthm, amsfonts, amscd}
\usepackage[english]{babel}
\usepackage{graphicx}
\usepackage{caption}
\usepackage{subcaption}
\usepackage{color,soul}
\usepackage{booktabs,dcolumn}                     
\usepackage[figureposition=bottom]{caption}       
\usepackage[usenames,dvipsnames,svgnames,table]{xcolor}
\usepackage{authblk}
\usepackage[title]{appendix}
\usepackage{tikz}
\usepackage{mathalfa}

\usepackage[utf8]{inputenc}
\usepackage[T1]{fontenc}
\usepackage{amsmath}
\usepackage{amsfonts}
\usepackage{amssymb}
\usepackage[version=4]{mhchem}
\usepackage{stmaryrd}
\usepackage{graphicx}
\usepackage[export]{adjustbox}

\usepackage{amsmath,amssymb,amsthm}
\usepackage{graphicx}
\usepackage{listings}
\usepackage{xcolor}
\usepackage[utf8]{inputenc}
\usepackage[T1]{fontenc}
\usepackage{amsmath,amsthm,amssymb,amsfonts}
\usepackage{geometry}

\newtheorem{thm}{Theorem}[section]

\usepackage{amsmath,amssymb}
\usepackage{graphicx}

\usepackage[colorlinks=true, allcolors=blue]{hyperref}

\usepackage{titlesec} 
\titleformat{\paragraph}
  {\normalfont\normalsize\bfseries}{\theparagraph}{1em}{}
\titlespacing*{\paragraph}{0pt}{3.25ex plus 1ex minus .2ex}{1.5ex plus .2ex}

\newcommand{\subsubsubsection}[1]{\paragraph{#1}}

\title{Polyhedral reconstruction via Boundary Control method}
\author{Dimitra Kyriakopoulou$^{\MakeLowercase{a}}$}
\affil{$^a$Mathematics Research Center, Academy of Athens, Greece}

\date{\vspace*{8mm}\textit{Dedicated to  Professor Yaroslav Kurylev}\vspace*{8mm}}

\begin{document}
\maketitle

\begin{abstract}
We study uniqueness of an elliptic Riemannian polyhedron using the elliptic version for Boundary Control method, which we presented in \cite{paper1}.
We also present interface detection criteria for hyperbolic Riemannian manifolds through introduction of the waveguide notion, the four-wave mixing notion, etc.
\end{abstract}

\section{Introduction}

In continuation of the work \cite{KuKi} where a uniqueness argument via Boundary Control (BC) method for a hyperbolic Riemannian polyhedron was presented, we proceed to uniqueness argument of an elliptic Riemannian polyhedron,
by extending the elliptic version for Boundary Control (BC) method, which we presented in \cite{paper1}, 
i.e. extending the pseudodifferential cone and edge singularities algebra that were used there.

The category of manifolds with singularities, denoted as $\mathfrak{M}_k$, includes $\mathfrak{M}_0$ for $C^{\infty}$ manifolds and progresses to higher orders representing more complex singularities, such as conical or edge singularities. This hierarchy allows for the definition of $\mathfrak{M}_{k+1}$ manifolds iteratively from $\mathfrak{M}_k$. Manifolds within $\mathfrak{M}_k$ support a natural differential operators algebra, characterized by principal symbolic hierarchies and ellipticity concepts. Prior research, including \cite{Downl_[120]}, \cite{Downl_[122]}, and \cite{Downl_[131]}, illustrates that advancing from level $k$ to $k+1$ in manifold analysis entails leveraging the parameter-dependent calculus of the preceding level, enriched by index theory and additional insights, challenging the notion of a simple inductive approach from $k$ to $k+1$. This complexity necessitates a detailed examination beyond mere induction, as suggested by \cite{SchCr}.
Hence, we prove

\begin{thm} \label{riempoly}
The BC method can provide a uniqueness argument for an elliptic Riemannian polyhedron.
\end{thm}

Furthermore, for the reconstruction of hyperbolic Riemannian polyhedron problem, we use the reconstruction per compartment choice; hence we resort to papers solving the problem in smooth manifolds, and we add interface and vertex detection conditions in the case of different wave types and different media.
For interface detection we try to find a condition equivalent to Snell's law but for perpendicular interface incidence.

Hence,  we show that

\begin{thm} \label{ewid}
For electromagnetic waves in isotropic dielectric media\\
(i) The amplitude transmission and reflection coefficients constitute conditions for interface detection.\\
(ii) Coupled mode equations constitute a vertex detection criterion.\\
(iii) For  polyhedron with multiple interfaces meeting at a vertex, the RC filter, i.e. a multi-directional coupler, can constitute a vertex detection criterion.
\end{thm} 

\begin{thm} \label{ewanl}
For electromagnetic waves in anisotropic non-linear media, e.g. a crystal, the four-wave-mixing constitutes a vertex detection criterion.
\end{thm}

\begin{thm} \label{acf}
For acoustic waves in fluid media \\
(i) the intensity transmission and reflection coefficients constitute conditions for interface detection.\\
(ii) The directional coupler with multiple transmission line matrices constitutes a vertex detection condition.
\end{thm}

\section{Reconstruction of an Elliptic Riemannian Polyhedron}

\noindent \textbf{Proof of Theorem \ref{riempoly}}.
Using the extension of Schulze's PDO edge algebra, which we used in \cite{paper1}, to higher singularities in the sense of Schulze again (see Appendix), we apply the method of \cite{paper1} to prove the Theorem. $\hfill \square$

\section[Interface and vertex detection conditions for Polyhedron]{Interface and vertex detection conditions for hyperbolic Riemannian Polyhedron reconstruction}

 We give conditions of interface and vertex detection for the reconstruction of a polyhedron with BC method in the case of different wave types and different media.
 
\subsection{Prerequisite facts}
\subsubsection{Riemannian polyhedron}

We define the polyhedron as in \cite{KuKi}. We start with a closed $n$-dimensional finite simplicial complex
\begin{equation*}
\mathcal{M}=\bigcup_{i=1}^I \mathrm{\Omega}^{(i)},   \label{(1)_3}
\end{equation*}
where $\mathrm{\Omega}^{(i)}$ are closed $n$-dimensional simpleces of $\mathcal{M}$. 
 
We assume that $\mathcal{M}$ is dimensionally homogeneous, i.e. any $k$-simplex, $0\leq k<n$, of $\mathcal{M}$ is contained in at least one $\mathrm{\Omega}^{(i)}$. We assume also that any $(n-1)$-dimensional simplex $\gamma$ belongs either to two different $n$ simpleces, $\mathrm{\Omega}^{(i)}$ and $\mathrm{\Omega}^{(j)}$, which, in this case, we often denote by $\mathrm{\Omega}_{-}$ and $\mathrm{\Omega}_+$, or to one $n$ simplex $\mathrm{\Omega}^{(i)}$. 
In the former case we call $\gamma$ an {\it interface} (sometimes $(n-1)$-dimensional interface) between $\mathrm{\Omega}_{-}$ and $\mathrm{\Omega}_+$, in the latter case we call $\gamma$ a {\it boundary} $(n-1)$-{\it simplex}, with $(n-1)$-simpleces having this property making the boundary $\partial \mathcal{M}$. 

We denote by $\mathcal{M}^{k}, 0\le k\le n$, the $k$-skeleton of $\mathcal{M}$ which consists of all $k$-simpleces of $\mathcal{M}$, with the differential structure on each $k$-simplex determined by its barycentric coordinates. Clearly, $\mathcal{M}=\mathcal{M}^{n}$. 

\subsubsection{Boundary Control (BC) method for hyperbolic problems}

The Boundary Control method (BCm), as elaborated e.g. in \cite{Bel}, offers a three-step process to construct a manifold that matches given data, involving defining input and state spaces (F and H, respectively), an input/state map \(W\), and an input/output map \(R\) that serves as the inverse data.

1. {\it Coordinatization}: Attach coordinates to each point in a manifold \(\Omega\), derived from state space \(H\), to create a metrically equivalent set \(\tilde{\Omega}\) based on reachable states and the mapping \(W\).
   
2. {\it Model Construction}: Use data \(R\) to determine inner products in the input space \(F\), leading to an auxiliary space \(\tilde{H}\) and mapping \(\tilde{W}\), which together form a model isometric to the original system's controllable part.

3. {\it Replication of \(\tilde{\Omega}\)}: By mimicking the coordinatization step with model states, replicate \(\tilde{\Omega}\), ensuring it is metrically and structurally identical to the original \(\Omega\).

Our metric is non-smooth between the polyhedron compartments. However, by smoothing $u$ with respect to time $t$, we can
extend it to non-smooth solutions, see e.g. \cite{[22]}.
 
 \subsection{Electromagnetic waves in isotropic dielectric media}
 The material of this section is from \cite{modopt}.
 
 \subsubsection{Interface detection}

We examine interactions of plane electromagnetic (EM) waves at the interfaces between two non-absorbing, isotropic, and homogeneous media. These media are characterized by their permittivity and permeability values denoted $(\varepsilon_{1},\mu_{1})$ and $(\varepsilon_{2},\mu_{2})$, for media 1 and media 2 respectively. 
Let $E_i$ and $H_i$, $i=1,2,3$ be the electric and magnetic fields associated with the incident, the refracted, and the reflected waves,  respectively.
 Let the incidence plane be the xz- plane.
 Let $\theta_1, \theta_2$ and $\theta 3$ be the angles of incidence, transmission/refraction, reflection, respectively, and $k_1, k_2$ and $k_3$ the
 propagation vectors along the incidence direction, the refracted direction and the
 reflected direction, respectively.  

Wave polarization at these interfaces is categorized into: P-polarized waves, where the magnetic field $\vec{H}$ is perpendicular to the incidence plane, also recognized as TM-polarized; and S-polarized waves, where the magnetic field lies within the incidence plane, hence termed TE-polarized. The amplitude transmission and reflection coefficients for both TM- and TE-polarized waves are defined as 
 \begin{equation} \label{t}
 t=\frac{E_2}{E_1},	
 \end{equation}
 and
 \begin{equation} \label{r}
 r=\frac{E_3}{E_1},	
 \end{equation}
indicating the wave amplitude ratios.

A lemma pertaining to normal incidence reveals that the amplitude transmission and reflection coefficients for TE- and TM-polarized waves are given, respectively, by 
 \begin{equation} \label{rprs}
 r_{p}=r_s=\frac{n_{1}-n_{2}}{n_{1}+n_{2}},
 \end{equation}
 and
 \begin{equation} \label{tpts}
 t_{p}=t_s=\displaystyle \frac{2n_{1}}{n_{2}+n_{1}},	
 \end{equation}
where $n_1$ and $n_2$ are the refractive indices of the respective media. This demonstrates that, at normal incidence, the coefficients are unaffected by wave polarization.\\

\noindent\textbf{Proof of Theorem \ref{ewid} (i).} Following a shortest geodesic, this will be normal to the interface from both of its sides. Hence we can recover waves up to a very small part past the interface. 
 We test the recovered values to spot the ones which satisfy the amplitude reflection and transmission coefficients given in the above Lemma, i.e. (\ref{tpts}) and (\ref{rprs}); these will correspond to the interface points. $\hfill \square$\\

\subsubsection{Vertex detection}
 
To detect the vertices we introduce waveguides, and make use of their coupling. For multi-dimensional polyhedron, we use directional couplers. Hence we start by introducing the relevant notions.
 
 \subsubsubsection{Waveguides}

Total internal reflection (TIR), a principle essential for understanding light behavior in optical fibers and other waveguide technologies, arises when light transitions from a denser medium to a rarer one ($n_1 > n_2$), resulting in complete reflection at the interface without transmission. This is captured for TE-polarized waves where the critical angle causes $\cos\theta_{2}$ to become imaginary, denoted by $\cos\theta_{2} = -i\sqrt{\frac{n_1^{2}\sin^{2}\theta_{1}}{n_{2}^{2}} - 1}$, indicating the absence of wave transmission across the boundary.

In TIR, the transmitted fields adjust to maintain $E_{y}, H_{x}, H_{z}$ components, with power flow exclusively parallel to the interface, and no cross-boundary energy transfer. This is mathematically expressed through wave equations for upward moving waves in the lower medium, $E_{1} = E_{0}e^{i(\omega t - \beta z - \kappa x)}$ and $E_{3} = E_{0}e^{i(\omega t - \beta z + \kappa x)}$, with $\beta$ and $\kappa$ describing wave propagation and boundary interaction.

Introducing two TIR interfaces creates a waveguide, confining the wave between them and allowing for evanescent decay outside, as shown by $\vec{E} = A\cos\kappa x e^{i(\omega t - \beta z)}$.

 \subsubsubsection{Waveguide Modes}

 Modes developed include TE Mode,  which is a TE-polarized wave guided in a waveguide, and the TM Mode, guided similarly as a TM-polarized wave. The waveguides can have different geometries -rectangular, circular, or arbitrary- with varied refractive index (RI) distribution $n(x, y)$. Specifically, for a planar design, RI depends on the $x$-coordinate, assuming the structure is infinitely extended along the $yz$-plane. Inhomogeneous wave equations for $E$ and $H$ fields in waveguides are 
 \begin{equation} \label{iw1}
 \nabla^{2}\vec{E}+\nabla(\frac{\nabla n^{2}}{n^{2}}\cdot\vec{E})=\mu_{0}\varepsilon_{0}n^{2}\frac{\partial^{2}\vec{E}}{\partial t^{2}},	
 \end{equation}
and
 \begin{equation} \label{iw2}
 \nabla^{2}\vec{H}+\frac{\nabla n^{2}}{n^{2}}\times(\nabla\times\vec{H})=\mu_{0}\varepsilon_{0}n^{2}\frac{\partial^{2}\vec{H}}{\partial t^{2}},
 \end{equation}	
respectively. These fields' variations depend on the waveguide's RI profile. TE and TM modes reveal specific field dependence of field components, with TEM mode equation for $E_y$ being 
\begin{equation} \label{TEM}
 \frac{\partial^{2}E_{y}}{\partial x^{2}}+(k_{0}^{2}n^{2}-\beta^{2})E_{y}=0,
 \end{equation}
illustrating mode behavior in RI variable waveguides. 

 \subsubsubsection{Mode coupling}

We examine the mode coupling phenomenon in axially invariant, uniform waveguides, where separate modes are supported with invariant field patterns during propagation. Unlike in uniform waveguides, where only the phase of modes changes, couplers induce variations in mode amplitudes through energy redistribution. This coupling may occur between modes within the waveguide or involve conversion from guided to radiation modes. Hence, we study these power redistributions in detail. Introducing a perturbation to an ideal waveguide facilitates the exchange of energy between modes, potentially enabling complete conversion from one mode to another under specific conditions. Mode coupling is a pivotal process that, under certain circumstances, allows for nearly total energy transfers between specific modes.

A directional coupler enhances mode exchange by facilitating co-directional coupling of two identical modes. It utilizes the waveguide's evanescent fields, which extend beyond the waveguide's boundaries. Thus, when two waveguides are placed sufficiently close that their evanescent fields overlap, energy is redistributed between the guides. With adequate distance, this fundamental power transfer can extend almost the entire interaction length, making power coupling and decoupling periodic functions of the interaction length.

The total field of the coupled waveguides, \(\psi(x,y,z)\), is expressible as a linear combination of the individual waveguide modes, \(\psi_1\) and \(\psi_2\), i.e., 
 \begin{equation} \label{p}
 \psi(x,y,z)=a(z)\psi_{1}(x,y)+b(z)\psi_{2}(x,y),	
 \end{equation}
The coupled mode equations are given by 
 \begin{equation} \label{ab1}
 \displaystyle \frac{\partial a}{\partial z}=-i(\beta_{1}+\kappa_{11})a-i\kappa_{12}b,
 \end{equation}
 \begin{equation} \label{ab2}
 \frac{\partial b}{\partial z}=
 -i(\beta_{2}+\kappa_{22})b-i\kappa_{21}a,	
 \end{equation} 
where \(\beta_1\) and \(\beta_2\) denote the propagation constants of the waveguides' modes, and the \(\kappa\) values represent the coupling coefficients.\\

\noindent\textbf{Proof of Theorem \ref{ewid} (ii).}  We consider two mathematical waveguides on two interfaces. The corner which constitutes the meeting point of the interfaces is going to satisfy the waveguide coupling mode equations; the smoothening of the metric allows for calculation of the values at the vertex points.. Let us note that we regard waveguides which have extremely small inter-interfacial area, hence coupled have again size approximately equal to the interface size.  $\hfill \square$\\
 
 \subsubsubsection{Multiple interfaces meeting at a vertex}
 The material of this section is in addition to \cite{modopt}, also from \cite{mult}.

In advancing the application of coupled mode equations to higher dimensions, it's essential to employ directional couplers, integral in enabling field interactions between dual-channel optical waveguides. These couplers facilitate energy transfer through evanescent mode coupling, leading to dynamic power redistribution between closely placed waveguides. The interaction is characterized by power tapping and reciprocal energy exchange, attributed to the overlapping of external fields via evanescent mode coupling. This process underpins the theoretical framework for analyzing power exchange between two waveguides, further elaborated through coupled mode equations.

Specifically, if we consider waves in two parallel waveguides with expressions 
 \begin{eqnarray}
 a(z)=a_{0}e^{-i\beta z}, \nonumber\\
 b(z)=b_{0}e^{-i\beta z}, \label{25.(2)}	
 \end{eqnarray}
 the resulting power in each waveguide can be respectively described as 
 \begin{equation}
 |a(z)|^2=1-\frac{\kappa^2}{\frac{1}{4}\Delta\beta^2+\kappa^2}\sin^2\{\left(\sqrt{\frac{1}{4}\Delta\beta^2+\kappa^2}\right)z\},	
 \end{equation}
 \begin{equation}
 |b(z)|^2=\frac{\kappa^2}{\frac{1}{4}\Delta\beta^2+\kappa^2}\sin^2\{\left(\sqrt{\frac{1}{4}\Delta\beta^2+\kappa^2}\right)z\}.	
 \end{equation}
 This mathematical formulation is crucial for understanding the variation of power along the z-axis in the waveguides.

Furthermore, an RC filter, constituted of multiple directional couplers and delay-line sections, acts as an add-drop filter. This arrangement can be mathematically represented through the multiplication of transmission matrices corresponding to each coupler and delay-line section. The transmission matrix for an N-stage RC filter is denoted as 
 \begin{equation} \label{(MP4.29)}
 \left[\begin{array}{l}
 \mathrm{Y}_{1}\\
 \mathrm{Y}_{2}
 \end{array}\right]=T_{c}(L_{N+1})\cdots T_{MZ}T_{c}(L_{2})T_{MZ}T_{c}(L_{1}) \left[\begin{array}{l}
 X_{1}\\
 X_{2}
 \end{array}\right],   	
 \end{equation}with $T_{MZ}$ and $T_c$ representing the transmission matrices for the delay-line section and the coupler, respectively.\\
 
\noindent\textbf{Proof of Theorem \ref{ewid} (iii).} So many stages should be added to the directional coupler as the number of interfaces that meet at a vertex minus 1. Then the directional coupler will perform coupling of the waveguides corresponding to these interfaces, and when this value is recovered this will indicate the vertex is detected. 
 $\hfill \square$\\
 
  \subsection{Electromagnetic waves in anisotropic non-linear media}
  The material of this section is from \cite{nonlinopt}.
  
 \subsubsection{Crystal detection}
 For the crystal lattice we need to detect the vertices, i.e. interface is not relevant. We use the four-wave mixing formula as the condition for checking for a vertex. First we use z-scan for third order nonlinearity susceptibility, which is necessary for the four-wave mixing formula.
 
 \subsubsubsection{Four-wave-mixing}

In four-wave mixing, three monochromatic pump fields of frequencies $\omega_1, \omega_2,$ and $\omega_3$ interact within a medium characterized by its third-order nonlinear susceptibility to produce a polarization $P_{NL}^{(3)}(\omega_s)$ at the mixed frequency $\omega_s = \omega_1 \pm \omega_2 \pm \omega_3$. This generates a signal field at $\omega_s$, propagating in the z-direction through a cubic or isotropic nonlinear optical medium. Under assumptions such as plane monochromatic fields, undepleted pump fields, and signal propagation along the z-direction, the signal field's evolution, based on slowly varying envelope approximation, is governed by the differential equation
\begin{equation}
\frac{dE_s}{dz} = -\frac{i\omega_s^2\mu_0\epsilon_0}{2k_s}\chi_{eff}^{(3)}E_1(z)E_2(z)E_3(z)e^{-i(\Delta\bar{k} \cdot \hat{s}_z)z},
\end{equation}
where  $\chi_{eff}^{(3)}$ is the third order nonlinear susceptibility and $\Delta\vec{k}=\vec{k}_s-\vec{k}_1-\vec{k}_2-\vec{k}_3$ represents the phase mismatch. Integration of this equation, given the initial condition $E_s(z=0)=0$, leads to the signal field expression:
\begin{equation}
E_s(z)=\frac{\omega_s^2\mu_0\epsilon_0}{2k_s(\Delta\vec{k} \cdot \hat{e}_z)}\chi_{eff}^{(3)}E_1(0)E_2(0)E_3(0)\left[\frac{\sin\left(\frac{(\Delta\vec{k} \cdot \hat{e}_z)z}{2}\right)}{\left(\frac{(\Delta\vec{k} \cdot \hat{e}_z)z}{2}\right)}\right]^2,
\end{equation}
highlighting that maximum signal intensity is achieved when the phase matching condition $\Delta k=0$ is satisfied. In a special case where the signal and one of the pump fields are degenerate ($\omega_s=\omega_3$), the signal field grows exponentially with distance, characterized by gain or loss, as described by
\begin{equation} \label{(NNO19.12)}
E_s(z)=E_s(0)\exp(g_sz),
\end{equation}
where $g_s(z)$ represents the gain or loss coefficient, emphasizing the controlled amplification or attenuation of the signal field through phase management in the nonlinear medium. \\
 
\noindent\textbf{Proof of Theorem \ref{ewanl}.} We consider the case $\omega_s=\omega_3$ above, where $E_3$ represents the normal geodesic on the vertex. We use Z-scan to get the third order nonlinear susceptibility, which is necessary for the four-wave mixing equation (\ref{(NNO19.12)}). The points on the normal geodesic that satisfy (\ref{(NNO19.12)}) indicate a vertex. $\hfill \square$
  
 \subsection{Acoustic waves in fluid media}
 The material of this section is from \cite{soupro}.
 
 The conditions for the interface detection use the pressure transmission and reflection
 coefficients. For the vertices, we use waveguiding, and an acoustic transmission line direction coupler.
 
 \subsubsection{Interface detection}

Regarding the propagation of acoustic waves across fluidic media interfaces, the process involves calculating the transmission and reflection coefficients for pressure waves. This calculation is crucial for understanding acoustic transmission lines or waveguides in various media. When an acoustic wave transmits from one medium to another, a portion of its energy is transmitted through the new medium while the remainder is reflected. This principle applies regardless of the medium being solid or fluid. In the case of lossless fluid media and planar waves, the complexity of this phenomenon is reduced. The fundamental equations for pressure \(p\) and velocity \(u\) in acoustic wave propagation are represented as 
 \begin{equation} \label{transm}
 \left\{\begin{array}{l}
 p(x,t)\\
 u(x,t)
 \end{array}\right\}=Re\left[ \left\{\begin{array}{ll}
 p_{+}(s) & p_{-}(s)\\
 \frac{p_{+}(s)}{z_{0}} & \frac{-p_{-}(s)}{z_{0}}
 \end{array}\right\}\left\{\begin{array}{l}
 e^{-sx/c}\\
 e^{sx/c}
 \end{array}\right\}\mathrm{e}^{st}\right],	
 \end{equation}
where \(z_{0}\) denotes the characteristic impedance. For a planar wave encountering a fluidic interface, the intensity transmission and reflection coefficients, \(T_{I}\) and \(R_{I}\), are used to quantify the energy distribution process. These coefficients are defined as 
 \begin{equation} \label{Eq.23.4b}
 T_{I}=\displaystyle \frac{4(Z_{2}/Z_{1})}{[(Z_{2}/Z_{1})-1]^{2}}   
 \end{equation}
 and
 \begin{equation} \label{Eq.23.5b}
 R_{I}=[\frac{(Z_{2}/Z_{1})-1}{(Z_{2}/Z_{1})+1}]^{2},  
 \end{equation} 
highlighting the role of characteristic impedances \(Z_{1}\) and \(Z_{2}\) in determining the wave's behavior at the interface. The boundary conditions, including equal acoustic pressure and the continuity of particle velocities across the interface, ensure material coherence and facilitate the precise analysis of wave interactions at media boundaries. This framework simplifies the complex phenomena associated with acoustic wave transmission and reflection, aiding in the understanding of acoustic wave behavior in different media.\\

\noindent\textbf{Proof of Theorem \ref{acf} (i).} Following a shortest geodesic, this will be normal to the interface from
 both of its sides. Hence we can recover waves up to a very small part past
 the interface. We test the recovered values to spot the ones which satisfy the transmission and reflection coefficients, 
 i.e. (\ref{Eq.23.4b}) and (\ref{Eq.23.5b}) in the above Lemma; these will correspond to the interface points. $\hfill \square$
 
 \subsubsection{Vertex detection}
 \textbf{Proof of Theorem \ref{acf} (ii).}
 The acoustic multi-stage directional coupler with as many transmission matrices as the number of interfaces meeting at a vertex minus 1, obeys the same equation (\ref{(MP4.29)}) as the one for the electromagnetic case, with the acoustic transmission line equations (\ref{transm}) replacing the electromagnetic ones. $\hfill \square$
 
\section*{Appendix: Algebra of higher singularities in the sense of Schulze}
\textbf{\textit{The iterative construction of higher singularities.}}
\cite{SchCr}, Section 10.3.1. 

In their work, Calvo, Martin, and Schulze \cite{Downl_[15]} introduce the category $\mathfrak{M}_k$ for singularity-ordered spaces. A space $M$ is part of this category if it progressively subtracts submanifolds $Y$, starting from a non-singular manifold ($\mathfrak{M}_0$), through to manifolds with singularities of order $k-1$. Specifically, for a space $M$ to belong to $\mathfrak{M}_k$, it must, after excluding a submanifold $Y$, fall within $\mathfrak{M}_{k-1}$. Furthermore, the dimensions of $Y$ indicate the nature of the singularity—zero dimensions imply a corner, while higher dimensions suggest an edge.

By deductive method, in [3] is established a hierarchy of $C^{\infty}$ submanifolds $Y^{(l)}$ within $M$, leading to a decomposition of $M$ into a sequence where each $M^{(j)}$ comprises submanifolds of decreasing singularity order, forming smooth edges of varying dimensions within $M$. 

Adequate differential operators of order $\mu$ for spaces within $\mathfrak{M}_k$ are denoted by $A\in \mathrm{D}\mathrm{i}\mathrm{f}\mathrm{f}_{\deg}^{\mu}(M)$, and are uniquely characterized by their behavior on non-singular portions of $M$ and near singular edges, with the latter described by a specialized formula involving the radial derivative and smooth coefficients. 

The concept of a principal symbolic hierarchy $\sigma(A)$ for these operators is introduced, capturing the traditional homogeneous principal symbol for non-singular parts of $M$ and extending it across the singularity order spectrum to describe the behavior of $A$ in relation to singular edges, culminating in a family of operators defined for the conical model space.\\

\noindent\textbf{\textit{Higher generations of weighted corner spaces.}}
\cite{SchCr}, Section 10.5.1. 

In the study of manifolds with complex singularities, the focus is on defining and understanding weighted Sobolev spaces, which are essential for modeling the manifold's geometric and analytic properties. These spaces, denoted as $\mathcal{K}^{s,\gamma}(X^{\wedge})$ for compact manifolds $X$ within a set $\mathfrak{M}_{k}$, and $\mathcal{W}^{s,\gamma}(X^{\wedge}\times \mathbb{R}^{q})$, serve to encapsulate the behaviors near singularities through the incorporation of weight tuples $\gamma \in \mathbb{R}^{k}$.

Manifolds $M$ are connected through a hierarchy of subspaces $M^{(j)}$, transitioning from complex structures to more simplified ones, ultimately allowing for an analytical framework that accommodates the manifold's singular characteristics. The construction of weighted Sobolev spaces is iteratively defined, starting from standard Sobolev spaces for $k=0$ scenarios and extending to weighted cone and edge spaces for manifolds with conical singularities or smooth edges, respectively. These definitions rely on stretching the manifold to fit the model spaces, with $\mathcal{H}^{s,\gamma}(\mathbb{M})$ and $\mathcal{W}^{s,\gamma}(\mathbb{M})$ being particular instances.

A remarkable aspect of these spaces is their invariance under natural chart transformations, facilitating the application across various manifold settings. Moreover, the transition to analyzing manifolds by directly addressing their singular nature, rather than their stretched counterparts, simplifies the notation and theoretical underpinnings, as illustrated by replacing $\mathcal{H}^{s,\gamma}(\mathbb{M})$ with $\mathcal{H}^{s,\gamma}(M)$.

The iterative process is further elucidated through the use of group isomorphisms, 
which underpin the local modeling near singular points and the seamless connection between local and global perspectives on the manifold. This leads to a comprehensive framework that not only addresses the manifold's inherent complexity but also provides a methodological basis for exploring weighted corner spaces in higher-dimensional settings.

Lastly, the theory accommodates both compact and non-compact manifolds, allowing for a broad application spectrum. This includes the consideration of manifolds as countable unions of compact sets and the introduction of spaces like $\mathcal{H}_{(\mathrm{c}\mathrm{o}\mathrm{m}\mathrm{p})}^{s,\gamma}(M)$ and $\mathcal{H}_{(\mathrm{lo}\mathrm{c})}^{s,\gamma}(M)$, which are tailored to handle the nuances of manifold structures with varying degrees of compactness. \\

\noindent\textbf{\textit{Additional edge conditions in higher corner algebras.}} \cite{SchCr}, Section 10.5.2. 

As we saw in \cite{paper1}
to extend an elliptic operator $A$ into a Fredholm operator across Sobolev spaces, formulating additional boundary conditions is essential, particularly by enhancing the boundary symbol into a set of isomorphisms. This may necessitate introducing vector bundles $J_{\pm}$ on the boundary, even for scalar operators, and applying to operators between distributional sections of vector bundles $E$ and $F$ to achieve Fredholm operators.

For manifolds with edges, an edge-degenerate operator $A$ necessitates supplementing the principal edge symbol $\sigma_{\wedge}(A)$ to transform it into a $2\times 2$ block matrix of isomorphisms for suitable vector bundles $J_{\pm}$. This process involves detailed constructions based on the bundles $E_{y}$ and $F_{y}$ derived from projecting vector bundles over the manifold's singularities and adjusting for admissible weights, which influence the choice of $J_{\pm}$ bundles.

Operators are considered Fredholm when they meet ellipticity conditions related to the edge algebra, signified by the bijectivity of the augmented edge symbol for all off-zero cotangent vectors. The Fredholm operators are then represented as $2\times 2$ block matrices connecting Sobolev spaces associated with the manifold and its boundary, based on specific weight data and the bundles $E$, $F$, $J_{-}$, and $J_{+}$.

For manifolds featuring hierarchical structure, the focus shifts to weighted Sobolev spaces associated with each subspace, characterized by their respective weights and vector bundles. The formulation extends to higher corner operator spaces, which encompass operators of a specific order and are represented as block matrices acting on the amalgamated weighted Sobolev spaces. These operators' ellipticity, crucial for their classification as Fredholm operators, is affirmed through the hierarchical structure of their principal symbols, accommodating both the manifold's interior and its corners.

\end{document}